\documentclass[a4paper,12pt]{article}
\usepackage[latin1]{inputenc}
\usepackage{amsmath,amssymb,amsthm,cite,verbatim}
\usepackage[all]{xy}
\usepackage[bookmarks=true]{hyperref}
\usepackage{color}
\usepackage{graphicx}

\title{Constraints on  the automorphism group of a  curve}
\author{ \ Josep Gonz\'alez  \footnote{The  author is partially supported by DGI grant  MTM2012-34611.
}}

\newtheorem{prop}{Proposition}%[section]
\newtheorem{lema}{Lemma}%[section]
\newtheorem{teo}{Theorem}%[section]
\newtheorem{cor}{Corollary}%[section]
\newtheorem{rem}{Remark}%[section]
\newtheorem{question}{Question}%[section]

\theoremstyle{definition}

\theoremstyle{remark}

\numberwithin{equation}{section}

\newcommand{\Q}{\mathbb{Q}}

\newcommand{\Z}{\mathbb{Z}}

\newcommand{\F}{\mathbb{F}}

\newcommand{\C}{\mathbb{C}}

\newcommand{\Gal}{\mathrm{Gal}}
\newcommand{\GL}{\operatorname{GL}}

\newcommand{\End}{\operatorname{End}}
\newcommand{\Aut}{\operatorname{Aut}}

\newcommand{\Frob}{\operatorname{Frob}}
\newcommand{\Jac}{\operatorname{Jac}}

\newcommand{\New}{\operatorname{New}}

\newcommand{\cS}{{\mathcal S}}

\newcommand{\cR}{{\mathcal R}}
\newcommand{\cQ}{{\mathcal Q}}

\date{}
\begin{document}
\maketitle

\begin{abstract}
For a curve  of genus $>1$ defined over a finite field, we present a sufficient criterion for  the non-existence of automorphisms of order a power of a rational prime. We show  how  this criterion can be used to determine the automorphism group of some modular curves of high genus.
\end{abstract}

\section{Introduction}
For a curve $X$ of genus $g>1$ defined over a field $K$,  the automorphism group  $\Aut_K (X)$ is finite.   In characteristic zero, it is  well-known that one has the Hurwitz bound
$|\Aut_K(X)|\leq 84(g-1)$.  In \cite{St},  the inequality  $|\Aut_K(X)|< 16 g^4$ is seen to hold in positive characteristic unless $X$ is a Hermitian curve. In particular, this provides a bound for the order of any automorphism of the curve. Nevertheless,
there is not a general procedure to discard possible orders.

We are interested in the case in which $K$ is a number field. The reduction of the curve $X$ at a prime of $K$ of  good reduction is a curve $\widetilde X$ defined over a finite field $\F_q$. Although $\Aut_{\F_q} (\widetilde X)$ may  strictly contain $\Aut_K(X)$, any information allowing to discard  orders of  the  elements  in the group $\Aut_{\F_q} (\widetilde X)$ will be useful for our goal. Moreover, if necessary, we can change the  prime of $K$ of good reduction for $X$.

The main result of this work is in Section 2.  For a curve $X$ of genus $>1$ defined over a finite field $\F_q$, we fix an integer $s>1$ which is a power of a rational prime. In Theorem \ref{crit},   we present a criterion which, under a certain condition on the sequence $\{|X(\F_{q^n})|\}_{n\geq 1}$ depending on $s$,  ensures   the non-existence of elements in $\Aut_{\F_q}(X)$ of order $s$. Although this criterion is not a characterization  for  the non-existence of such  automorphisms, it is certainly a powerful tool that can be applied in many situations. In Section 3, in order to show  the  efficacy of this tool, we apply it  to determine the automorphism groups of  some modular curves.
%============================================================================================================
\section{Automorphisms of a curve defined over a finite field}
Let $\F_q$ be the finite field with $q$ elements and let $X$ be a curve  of genus $g>1$ defined over $\F_q$.
\begin{teo}\label{crit}  If for a rational prime $N$ and an integer $m>0$,  the sequence of integers $\{P_{N^m}(n)\}_{n\geq 1}$  defined by
$$
0\leq P_{N^m}(n)\leq N^m-1\,,\quad P_{N^m}(n):= \left\{\begin{array}{rr}|X(\F_q)|\pmod{N^m} & \text{if $n=1$,}\\[6 pt] |\cup_{i=1}^{n}X(\F_{q^i})|-|\cup_{i=1}^{n-1}X(\F_{q^i})| \pmod{N^m} &\text{if $n>1$,}\end{array}\right.
$$
satisfies the condition
\begin{equation}\label{criterion}\sum_{n\geq  1}P_{N^m}(n)>\displaystyle{\left\lfloor\frac{2\,g\, }{N-1}\right\rfloor+ 2\frac{N^m-1}{N-1}}\,,
\end{equation}
then $X$ does not have any automorphisms  defined over $\F_q$ of order $N^m$.
\end{teo}
\noindent {\bf Proof.}  Assume  that there exists $u\in\Aut_{\F_q}(X)$ of order $N^m$. Let $G$ be the subgroup of $\Aut_{\F_q}(X)$ generated by $u$. Let  $\cR$ denote the set of ramification points of the natural projection $\pi_G\colon X\rightarrow X/G$.

The group $G$ acts on the set $X(\F_{q ^{n}})$ as a subgroup of permutations. In particular, if  $Q\in X(\F_{q^{ n}})$, then the orbit of $Q$ under $G$, say $\cS_Q$,  is contained in   $ X(\F_{q^n})$. It is clear that $|\cS_Q|=N^m$  if, and only if, $Q\notin \cR$. Let $n_0>0$ be an integer such that $\cR\subset X(\F_{q^{n_0}})$. Hence,  we have
 $$ |X(\F_{q^{n_0}})|\equiv |\cR|\pmod{N^m}\,.$$
 So the sequence $A(n):=|\cup_{i=1}^nX(\F_{q^{ i}})|$, $n\geq 1$, satisfies
   the condition   $A(n)\equiv |\cR|\pmod {N^m}$ for all $n\geq n_0$. Therefore,
\begin{equation}\label{cR}
\sum_{n\geq  1}P_{N^m}(n)=\sum_{n\geq 1}^{n_0}P_{N^m}(n)\leq |\cR\cap X(\F_q)|+\sum_{n=2}^{n_0}|\cR\cap(\cup_{i=1}^{n}X(\F_{q^ i})-\cup_{i=1}^{n-1}X(\F_{q^ i}))|=|\cR|\,.
\end{equation}

For $0\leq i\leq m-1$, let  $\cR_i$ denote the subset of $\cR$ consisting of the points whose isotropy group in $G$ is the subgroup  generated by $u^{ N^i}$,  and set $r_i:=|\cR_i|$. One has  $\cR=\cup_{i=0}^{m-1}\cR_i$ and, moreover, the ramification index of $\pi_G$ at a point $Q\in\cR_i$  is the order of its isotropy group, i.e. $N^{m-i}$. By the Riemann-Hurwitz formula applied to $\pi_G$, we obtain
$$
 N^m(2g_G-2)+(N^m-1)r_0+\cdots +(N^i-1)r_{m-i}+\cdots +(N-1) r_{m-1}\leq 2 g-2\,,
$$
where $g_G$ is the genus of   $X/G$. In particular, we have
$$
 (N^m-1)r_0+\cdots +(N-1) r_{m-1}\leq 2g+2 (N^m-1) \,.
$$
Therefore,
  \begin{equation}\label{des}
 |\cR|= r_0+ \cdots + r_{m-1}\leq \frac{1}{N-1} \left((N^m-1)r_0+\cdots +(N-1) r_{m-1}\right)\leq\frac{ 2g+2 (N^m-1)}{N-1}\,.
 \end{equation}
Combining    the inequalities (\ref{cR}) and (\ref{des}), we get
\begin{equation}\label{nec}\sum_{n\geq  1}P_{N^m}(n)\leq \displaystyle{\left\lfloor\frac{2\,g\, }{N-1}\right\rfloor+ 2\frac{N^m-1}{N-1}} \,,
\end{equation} which proves the statement. \hfill $\Box$

\begin{rem}  To apply Theorem \ref{crit}, we only need to know the characteristic polynomial $Q(x)$ of $\Frob_q$ acting on the Tate module of $\Jac (X)$. Indeed, if $\alpha_1,\cdots, \alpha_{2g}$ are the roots of $Q(x)$, then
$$
|X(\F_{q^n}) |=1+ q^n-\sum_{i=1}^{2g} \alpha_i^n\,.
$$
For $n>1$, the integer $R(n):=|\cup_{i=1}^{n}X(\F_{q ^i})|-|\cup_{i=1}^{n-1} X(\F_{q^i})|$ can be computed from the sequence $\{|X(\F_{q^i})|\}_{1\leq   i\leq n}$  as follows. Let $\{p_1, \cdots ,p_k\}$ be  the set of   primes dividing $n$ and  put $d_i=n/p_i$  for $1\leq i\leq k$. Then,
\begin{equation}\label{easy}
R(n)=|X(\F_{q^{n}})|-\sum_{r=1}^k(-1)^{r+1}
\sum_{1\leq i_1<\cdots <i_r\leq k}|X(\F_{q^{\gcd(d_{i_1},\cdots,d_{i_r})}})|\,.
\end{equation}
Indeed,  using that
$$X(\F_{q ^{d_1}})\cap X(\F_{q ^{d_2}})=X(\F_{q ^{\gcd(d_1,d_2)}})\,,\quad\text{ and if }
d_1|d_2 \text{ then }X(\F_{q ^{d_1}})\cup X(\F_{q ^{d_2}})=X(\F_{q ^{d_2}})\,,
$$
we obtain
$$
\begin{array}{rl}
R(n)&= |X(\F_{q^{n}})|-|X(\F_{q^{n}})\cap \left(\cup_{i=1}^{n-1} X(\F_{q^i})\right)|\\[6 pt]
     &=|X(\F_{q^{n}})|-|\cup_{d|n, d<n} X(\F_{q^d})|=|X(\F_{q^{n}})|-|\cup_{i=1}^k X(\F_{q^{d_i}})|\\[6 pt]
   &=|X(\F_{q^{n}})|-\sum_{r=1}^k(-1)^{r+1}
\sum_{1\leq i_1<\cdots <i_r\leq k}|X(\F_{q^{\gcd(d_{i_1},\cdots,d_{i_r})}})|\,.
\end{array}
$$
Note that, if $\ell_1=2< \cdots <\ell_r$ are the first $r$ rational primes, for $n<\prod_{i=1}^r \ell_i$,
the sum  given in  (\ref{easy})  contains at most $2^{r-1}$ terms.

To be more precise, to apply Criterion \ref{crit} we only need to know $Q(x)\pmod{N^m}$. In other words, we can change the polynomial $Q(X)$ by a polynomial $T(x)\in\Z[x]$ such that $Q(x)\equiv T(x)\pmod {N^m}$. We can      determine $R(n)\pmod{N^m}$ from the roots of $T(x)$ by applying  the  procedure described for the roots of $Q(x)$.
\end{rem}

\begin{rem}\label{r} If we can prove that, for an automorphism $u\in\Aut_{\F_q}(X)$ of order $N^m$, there exists an integer $r$ such that $|\cR| \leq r< \left\lfloor\frac{2\,g\, }{N-1}\right\rfloor+ 2\frac{N^m-1}{N-1}$, then  the condition (\ref{criterion}) in Theorem \ref{crit} can be replaced with the condition
$$
\sum_{n \geq 1}P_{N^m}(n)> r\,.
$$
\end{rem}
\begin{rem}
  The non-existence of an automorphism in $\Aut_{\F_q}(X)$ of order $N^m$ is not a necessary   condition  to satisfy   condition (\ref{criterion}). For instance, it may be that two non-isomorphic curves $X$ and $Y$ defined over $\F_q$ have jacobians which are isogenous over $\F_q$. If   $Y$ has an automorphism of order  $N^m$, then  condition (\ref{criterion}) is not satisfied, even if $X$ does not have an automorphism of order $N^m$. Also, if the group $G =\Aut_{\F_q}(X)$ is nontrivial, then   one has $ |\cup_{i=1}^n X(\F_{q^i})|\equiv |\cS|\pmod {|G|}$ for almost all $n$, where
  $\cS$ is the set of  ramification points of the covering $X\rightarrow X/G$.  It may be that condition (\ref{criterion}) is not satisfied   when we take  $N^m$ dividing $|G|$ and $G$ does not contain any $N^m$-cyclic subgroup.
\end{rem}
\begin{question} For a prime $N$, the condition that the sequence
$ \{|\cup_{i=1}^n X(\F_{q^i})|\pmod {N^m} \}_{n\geq n_0}$ is constant for some integer $n_0$ seems to be   strong.  If there exists a curve $Y$ defined over $\F_q$ such that $\Jac (Y)$ and $\Jac (X)$ are isogenous over $\F_q$ and the order of the group $\Aut_{\F_q}(Y)$ is a multiple of $N^m$, then this condition is satisfied. Is the converse true?
\end{question}

Several consequences can be obtained from Theorem \ref{crit}. Next, we present two of them.
\begin{cor}   If there is an integer $n_0>0$ such that
$$  2< |\cup_{i=1}^{n_0+1}X(\F_{q^i})|-|\cup _{i=1}^{n_0}X(\F_{q^i})| <2g+2\,,$$
then there are not any automorphisms in $\Aut_{\F_q}(X)$ of order a  prime $N>2g+1$, which improves  the result obtained through the Hurwitz bound.
\end{cor}
\begin{cor} If $u\in\Aut _{\F_q}(X)$ has order $N^m$, then $\sum
_{n \geq 1} P_{N^m}(n)$ is a lower bound for the cardinality of the set  of ramified points of the covering $X\rightarrow X/G$, where $G$ is the subgroup of $\Aut _{\F_q}(X)$ generated by $u$.

\end{cor}

\section{Application to some modular curves}
   Let $\New_N$ denote   the set of normalized newforms in   $S_2(\Gamma_0(N))^{\operatorname{new}}$ and let $\New_N^+$ be the set  $\{f\in \New_N\colon w_N(f)=f\}$, where $w_N$ is  the Fricke  involution.
 For $f\in\New_N$, let $S_2(f)$ be the $\C$-vector space  of cusp forms spanned  by  $f$ and its  Galois conjugates. Let us denote by   $A_f$ the abelian variety attached to $f$ by Shimura. It is  a quotient of $J_0(N):=\Jac (X_0(N))$  defined over  $\Q$ and the  pull-back of $\Omega^1_{A_f/\Q}$  is the $\Q$-vector subspace of elements in  $S_2(f) dq/q$ with rational $q$-expansion, i.e. $S_2(f) dq/q\cap\Q[[q]]$, where $q=e^{2\pi\,i\,z}$ for $z$ in the complex upper half-plane.
 Moreover, the endomorphism algebra $\End_\Q^0(A_f):=\End_\Q(A_f)\otimes \Q$ is isomorphic to  a totally real number field $E_f$ whose degree is equal to $\dim A_f$.

Let $G_\Q$ denote the absolute Galois group $\Gal (\overline{\Q}/\Q)$.  Let $X$ be a curve  of genus $g>0$ defined over $\Q$ such that $\Jac (X)$  is a quotient of the jacobian of the curve $X_0(N)$ defined over $\Q$.  There exists a subset $\cS$ of the set $\cup_{M|N} \New_M$, which is stable under Galois conjugation,  such that  $\Jac (X)$ is isogenous over $\Q$ to the  abelian  variety $\prod_{f\in \cS/G_{\Q}}A_f^{n_f}$ for some integers $n_f>0$.  If $\ell$ is a prime of good reduction for $X$ not dividing $N$, by the Eichler-Shimura congruence, we can compute the characteristic polynomial $Q(x)$ of  $\Frob_\ell$ acting on the Tate module of $\Jac (X\otimes \F_\ell)$ through the  $\ell$-Fourier coefficients $a_\ell(f)$ of  the newforms $f$ in $\cS$:
$$
Q(x)=\prod_ {f\in\cS}(x^2-a_\ell(f) x+\ell)^{n_f}\,.
$$

\subsection{The split Cartan modular curves $X_{s}(p)$} For a prime $p$, let us denote by $X_s(p)$ the modular curve attached to the normalizer of a split Cartan subgroup of $\GL_2(\F_p)$. This curve is a quotient of the modular curve $X(p)$ defined over $\Q$ and  is isomorphic over $\Q$ to the modular curve $X_0^+(p^2)=X_0(p^2)/\langle w_{p^2}\rangle$. In \cite{go15}, the automorphism group of the curve  $X_{s}(p)$  is determined  for all primes $p$.
There, to conclude the article, it is needed to prove that $X_{s}(p)$ does not have any involutions defined over the quadratic field $K=\Q(\sqrt{p^*})$, where $p^*=(-1)^{(p-1)/2}$, for $p=17,19,23,29,31$. In fact, this problem is at the origin of Theorem 1 for $N^m=2$. Since in \cite{go15} it is proved that the number of fixed points of an involution  is $\leq 12$, we can apply Remark \ref{r} to the reduction of $X_s(p)$ at a prime of $K$ over a rational prime $\ell\neq p$, i.e.  we can use that the condition $\sum _{n}P_2(n)>12$  implies the non-existence of an involution of $X_s(p)$ defined over $K$. This fact is proved by taking  a prime $\ell$ splitting in $K$ and  by applying  this  version of Theorem \ref{crit}  to the curve $X_{s}(p)\otimes \F_\ell$ for $N^m=2$ (see Section 5 of \cite{go15}).
%$$
%\begin{array}{|c|r|c|}
%p& g&  \sum_nP_2(n)\\\hline
%17 &   &\sum_{n\leq 52}P_2(n)=15\\[6 pt]
%19 &   &\sum_{n\leq 52}P_2(n)=15\\[6 pt]
%23 &   &\sum_{n\leq 52}P_2(n)=15\\[6 pt]
%29 &   &\sum_{n\leq 52}P_2(n)=15\\[6 pt]
%31 &   &\sum_{n\leq 52}P_2(n)=15
%\end{array}
%$$

\subsection{The modular curves $X_0^+(p)$} In \cite{BH},  the automorphism group of the modular curves $X_0^+(p):=X_0(p)/\langle w_p\rangle$  is determined for all primes $p$.  After applying some theoretical results and  to conclude the article, the authors  need to prove that the modular curve $X_0^+(p)$ does not have any involution defined over $\Q$ for $p=163,193,197,211,223,227,229,269,331,347,359,383, 389,431,461,563,571,607$. In order to do that, they apply two different arguments. The first one  is used to discard $11$ cases and the second one allows to discard the remaining $7$ cases.  Although in \cite{BH} it is proved that the number of fixed points of an involution is $\leq 12$, next we show the table obtained by applying Theorem \ref{crit} (without using   Remark \ref{r}) to the curve $X=X_0^+(p)\otimes\F_2$ and $N^m=2$:
$$
\begin{array}{|c|r|c||c|r|c||c|r|c|}
p& g&  \sum_nP_2(n) & p& g&  \sum_nP_2(n)&p& g&  \sum_nP_2(n)\\ \hline
163 & 6 &\sum_{n\leq 53}P_2(n)=15 &229 & 7 &\sum_{n\leq 63}P_2(n)= 17 &389 &11 &\sum_{n\leq 123 }P_2(n)=25 \\[6pt]
193 & 7 &\sum_{n\leq 58}P_2(n)=17 &269 & 6 &\sum_{n\leq 43}P_2(n)= 13 & 431 & 8 &\sum_{n\leq 89}P_2(n)=19 \\[6pt]
197 & 6 &\sum_{n\leq 42}P_2(n)=15 &331 & 11 &\sum_{n\leq 79}P_2(n)= 25 & 461 & 12 &\sum_{n\leq 99}P_2(n)=27\\[6pt]
211& 6 &\sum_{n\leq 60}P_2(n)=15 &347 & 10 &\sum_{n\leq 74}P_2(n)=23  &563 & 15 &\sum_{n\leq 116}P_2(n)= 33\\[6pt]
223 & 6 &\sum_{n\leq 54}P_2(n)=15 &359 & 6 &\sum_{n\leq 60}P_2(n)= 15  & 571&  19&\sum_{n\leq 156}P_2(n)=41\\[6pt]
227 & 5 &\sum_{n\leq 40}P_2(n)= 13&383 & 8 &\sum_{n\leq 88}P_2(n)=19  & 607 & 19 &\sum_{n\leq 166}P_2(n)= 41
\end{array}
$$
In all cases $\sum_{n}P_2(n)>2\, g+2$ and, thus, all these curves do not  have any involutions defined over $\Q$.

\subsection{The non-split Cartan modular curves $X_{ns}(p)$}

Let $p$ be a rational prime and let  $X_{ns}(p)$ be the  modular curve attached to a non-split Cartan  subgroup of $\GL_2(\F_p)$. This curve is a quotient of the modular curve $X(p)$ defined over $\Q$, which has a canonical involution $w$ defined over $\Q$, the so-called modular involution. The genus $g$ of $X_{ns}(p)$ is greater than $1$ for $p\geq 11$. In \cite{DFGS}, the following  is proved
$$\Aut(X_{ns}(11))=\Aut_\Q(X_{ns}(11))\simeq (\Z/2\Z)^2\,.$$
 In \cite{do15}, it is proved that for  $p\geq 37$  all automorphisms of $X_{ns}(p)$ preserve cusps and, moreover, if $p\equiv 1\pmod{12}$ then $\Aut(X_{ns}(p))=\{1,w\}$.

It is expected that $\Aut(X_{ns}(p))=\{1,w\}$ for $p>11$.  The goal of this subsection is to prove this fact
 for $ 13\leq p \leq 31$. We point out that the genera of these six curves are $8, 15, 20, 35, 54$ and $63$.

\vskip 0.2 cm
Set $X_{ns}^+(p)=X_{ns}(p)/\langle w\rangle$ and let us denote by $g^+$ its genus. For $p\geq 11$,  the splitting  over $\Q$ of the jacobians of these curves is as follows (cf. \cite{Chen}):
$$
J_{ns}(p):=\Jac(X_{ns}(p))\stackrel{\Q}\sim \prod_{f\in \New_{p^2}/G_\Q} A_f\,,\quad J_{ns}^+(p):=\Jac(X_{ns}^+(p))\stackrel{\Q}\sim \prod_{f\in \New_{p^2}^+/G_\Q} A_f\,.
$$

From now on,
   $\chi$  denotes the  quadratic Dirichlet character of conductor $p$, i.e. the Dirichlet character  attached to the quadratic number field  $K=\Q(\sqrt{p^*})$, where $p^*=(-1)^{(p-1)/2}$. Next, we summarize some facts concerning  the modular abelian varieties $A_f$ attached to  newforms $f\in\New_{p^2}$ (see Section 2 of  \cite{go15} for  detailed references).

 The map $f\mapsto f\otimes \chi$ is a permutation of the set $\New_{p^2}\cup\New_p$. Under this bijection, there is a unique newform $f$, up to Galois conjugation, such that $f=f\otimes \chi$ when $p\equiv 3 \pmod 4$ and, moreover, in this case $f\in\New_{p^2}$.

If $f\in\New_{p^2}$ has complex multiplication (CM), i.e. $f=f\otimes\chi$, then the dimension of $A_f$ is the class   number of $K$ and  $A_f$ has all its endomorphisms defined over the Hilbert class field of $K$. The endomorphism algebra  $\End_K^0(A_f)$ is isomorphic to the CM field $E_f\otimes K$ which only contains the roots of unity $\pm 1$. Moreover, $f\in\New_{p^2}^+$ if, and only if, $p\equiv 3\pmod 8$.

 Let $f=\sum a_n q^n  \in \New_{p^2}$ be without CM. If $f$ has an inner twist $\chi'\neq 1$, i.e. $f\otimes \chi'={}^{\sigma} f$ for some $\sigma\in G_\Q$, then $\chi'=\chi$ because $\chi'$ must be a quadratic character of conductor dividing $p^2$.
In such a  case,  $\End ^0(A_f)=\End_K^0(A_f) $ is a noncommutative algebra.
 More precisely, set $F_f:=\Q( \{ a_\ell^2\})$, with $\ell$ running over the set of all rational  primes. If $A_f$  is simple, then  $\End_K^0(A_f) $ is a  quaternion algebra $\cQ_f$ over $F_f$ (QM case), otherwise $A_f$ is isogenous over $K$ to the square of an abelian variety $B_f$ and $\End_K^0(A_f) $ is isomorphic to the matrix algebra $M_2(F_f)$ (RM case).

If $\chi$ is not an inner twist for $f\in\New_{p^2}$, then  $A_f$ is simple  and $\End^0 (A_f)$ is isomorphic to $E_f$ (RM case).

For two distinct $f_1, f_2\in \New_{p^2}/G_\Q$, the abelian varieties $A_{f_1}$ and $A_{f_2}$ are not isogenous over $\Q$ and are isogenous if, and only if, $f_1\otimes \chi={}^{\sigma} f_2$ for some $\sigma\in G_\Q$. In this particular case, there is an isogeny defined over $K$.

In the sequel, we restrict our attention to the values  $13\leq p\leq 31$. The next lemma can be obtained  through  the instruction {\bf BrauerClass} in the program {\it Magma}.
\begin{lema}
There is not any $f\in\New_{p^2}$ with   quaternionic multiplication for  all $13\leq p\leq 31$.
\end{lema}

\vskip 0.2 cm
Let us fix a set $\{f_1,\cdots,f_r\}$ of representative cusp forms for  the set $\New_{p^2}/G_\Q$. We introduce the  subsets $\cS_{cm}$,  $\cS_{rm}$, $\cS_{s}$ and $\cS_{t}$ of $\New_{p^2}/G_\Q$ as follows.
The subsets $\cS_{cm}$ and $\cS_{rm}$ are the sets of newforms in $\New_{p^2}/G_\Q$ having $\chi$ as an inner twist and  corresponding to the  CM and RM cases respectively. The subsets $\cS_s$ and $\cS_t$ are defined as follows:
$$ \cS_{s}=\{ f\in \New_{p^2}/G_\Q\colon f\otimes \chi \in\New_p/G_\Q\}\,,\,\,\cS_{t}=\{f_i \in \New_{p^2}/G_\Q\colon f_j=f_i\otimes \chi, i<j\}\,.
$$
For $\New_{p^2}^+/G_\Q$, we introduce the following four sets $\cS_{cm}^+=\cS_{cm}\cap \New_{p^2}^+$, $\cS_{rm}^+=\cS_{rm}\cap \New_{p^2}^+$, $\cS_s^+=\cS_s\cap \New_{p^2}^+$ and
$ \cS_{t}^+=\{f_i \in \cS_t\cap \New_{p^2}^+\colon f_i\otimes \chi\in \New_{p^2}^+\}$.  Hence,
the splitting over $K$ of $J_{ns}(p)$ and $J_{ns}^+(p)$  are
$$
J_{ns}(p)\stackrel{K}\sim \prod_{f\in \cS_{cm}} A_f\prod_{f\in \cS_{s}}A_f\prod_{f\in \cS_{rm}} B_f^2\prod_{f\in\cS_t} A_f^2$$
 and
 $$J_{ns}^+(p)\stackrel{K}\sim \prod_{f\in \cS_{cm}^+} A_f\prod_{f\in \cS_{s}^+}A_f\prod_{f\in \cS_{rm}^+} B_f^2\prod_{f\in\cS_t^+} A_f^2\,.
$$
The corresponding decomposition of their endomorphism algebras over $K$  are
\begin{equation}\label{dc}
\End_K^0(J_{ns}(p))\simeq \prod_{f\in \cS_{cm}} E_f\otimes K\prod_{f\in \cS_{s}}E_f\prod_{f\in \cS_{rm}} M_2(F_f)\prod_{f\in\cS_t} M_2(E_f)
\end{equation}
and
\begin{equation}\label{dc+}
\End_K^0(J_{ns}^+(p))\simeq \prod_{f\in \cS_{cm}^+} E_f\otimes  K\prod_{f\in \cS_{s}^+}E_f\prod_{f\in \cS_{rm}^+} M_2(F_f)\prod_{f\in\cS_t^+} M_2(E_f)\,.
\end{equation}
For $13\leq p\leq 31$, the following table shows the description of the sets $\New_{p^2}/G_\Q$ and $\New_{p^2}^+/G_\Q$ as well as the action of the map $f\mapsto f\otimes \chi$ on the set $(\New_{p^2}\cup\New_p)/G_\Q$.

$$
\begin{array}{c|c|c|c|c|c|c|c|c|}
p & g& \New_{p^2}/G_\Q & \dim A_{f_i}& \cS_{rm}& \cS_{cm} & \begin{array}{c}\cS_t\\ \text{(twists)}\end{array} &g^+& \New^+_{p^2}/G_\Q\\\hline\hline
13 &8 &\{f_1,f_2,f_3\} &\left\{\begin{array}{cr}2\,, & i=1\\
  3 \,,&2\leq i\leq 3\end{array}\right. &\{ f_1\}& \emptyset &\begin{array}{c}\{ f_2\}\\f_3=f_2\otimes \chi\end{array} & 3 & \{ f_2\}\\\hline
17 &15  &\{f_1,\cdots,f_6\} & \left\{\begin{array}{cr}1\,, & i=1\\
  2\,,&2\leq i\leq 3\\
  3\,,& 4\leq i\leq 5\\
  4\,,& i=6\end{array}\right. & \{f_6\}&\emptyset &\begin{array}{c}\{f_2,f_4\}\\ f_3=f_2\otimes \chi \\
  f_5=f_4\otimes \chi\end{array}& 6 &\{f_1,f_2,f_4\}  \\\hline
19 & 20&\{f_1,\cdots,f_9\}& \left\{\begin{array}{cr} 1\,, & 1\leq i \leq  2\\
2\,, & 3\leq i \leq 6\\
3\,, &7\leq i \leq 8\\
4 \,,& i=9\end{array}\right.& \{f_9\}& \{f_1\}& \begin{array}{c}\{f_3,f_5,f_7\}\\ f_4=f_3\otimes \chi \\
  f_6=f_5\otimes \chi\\ f_8=f_7\otimes \chi\end{array}& 8& \{f_1,f_7,f_9\} \\\hline
23 &31 &\{f_1,\cdots,f_{10} \}&
\left\{\begin{array}{cr} 2\,, & 1\leq i \leq  5\\
3\,, & i= 6\\
4\,, &7\leq i \leq 8\\
5 \,,&9\leq i\leq 10 \end{array}\right. &\{f_7,f_8\} & \{f_6\}&  \begin{array}{c} \{f_1,f_4,f_9\}\\f_2=f_1\otimes \chi \\
  f_5=f_4\otimes \chi\\ f_{10}=f_9\otimes \chi\end{array}& 13&\{ f_7,f_8,f_9 \}\\\hline
29 & 54&\{f_1,\cdots,f_{11} \}&\left\{\begin{array}{rr} 2\,, & 1\leq i \leq  4\\
3\,, & 5\leq i \leq 6\\
6\,, &7 \leq i \leq 8\\
8 \,,& 9\leq i \leq 10\\
12\,, & i=11
\end{array}\right. & \{f_3,f_{11}\}&\emptyset  &\begin{array}{c}\{f_1,f_5,f_7,f_9\}\\f_4=f_1\otimes \chi\\f_6=f_5\otimes\chi\\f_8=f_7\otimes \chi\\f_{10}=f_9\otimes \chi\end{array} &24 &\begin{array}{c}\{f_1,f_2,f_5,\\f_6,f_7,f_9 \}\end{array} \\\hline
31 & 63&\{f_1,\cdots,f_{12}\} & \left\{\begin{array}{rr} 2\,, & 1\leq i \leq  6\\
3\,, & i =7\\
4\,, &  i = 8\\
 8\,,& 9\leq i \leq 10\\
12\,, & i=11\\
16\,,& i=12
\end{array}\right.&\begin{array}{c}\{f_5,f_8,\\f_{11},f_{12}\}\end{array} &\{f_7 \}&\begin{array}{c}\{f_1,f_2,f_9\}\\f_4=f_1\otimes \chi\\f_6=f_2\otimes\chi\\f_{10}=f_9\otimes \chi\end{array} & 28&\begin{array}{c}\{f_1,f_2,\\f_9,f_{12}\}\end{array}\\ \hline
\end{array}
$$
\begin{centerline} {Table 1}
\end{centerline}

\vskip 0.2 cm
\noindent
The label of the newforms in $\New_{p^2}$ is the one given by {\it Magma}. For a prime $p$, the set $\cS_s$ is the set of newforms $f$ which do not appear in the columns corresponding to $\cS_{rm}$, $\cS_{cm}$  and  $\cS_t$ (twists).

\vskip 0.3cm
\begin{prop} Let $p$ be a prime such that $13 \leq p\leq 31$. Then,
\begin{itemize}
\item[(i)] The group $\Aut(X_{ns}^+(p))$ is trivial.
\item[(ii)] The modular involution $w$ is the only nontrivial automorphism of $X_{ns}(p)$.
\end{itemize}
\end{prop}

\noindent{\bf Proof.} For $p=13$, we already know that $\Aut (X_{ns}^+(13))$ is trivial because $X_{ns}^+(13)$ is not hyperelliptic (cf. \cite{Baran}) and the endomorphism algebra $\End^0 (J_{ns}^+(13))$ is a totally real number field which only contains the roots of unity $\pm 1$.

We split the proof into the following steps.

\vskip 0.3 cm

\noindent {\it Step 1: All automorphisms of $X_{ns}(p)$ and $X_{ns}^+(p)$ are defined over $K$. }

On the one hand, for two distinct $f_1,f_2$ lying in  $\New_{p^2}/G_\Q$,  without CM, $A_{f_1}$ and $A_{f_2}$ are isogenous if, and only if, $f_2$ is a Galois conjugate of $f_1\otimes \chi$ and, in this case, the isogeny is defined over $K$.
On the other hand, if $f\in\New_{p^2}/G_\Q$   does not have CM, all endomorphisms of  $A_f$ are defined over  $K$. Hence, if  $\New_{p^2}/G_\Q$, resp. $\New_{p^2}^+/G_\Q$, does not contain a newform with CM, all endomorphisms of $J_{ns}(p)$, resp. $J_{ns}^+(p)$, are defined over $K$ and, in particular, also all automorphisms of the corresponding curve.

Assume that $\New_{p^2}/G_\Q$, resp $\New_{p^2}^+/G_\Q$, contains a newform $f$ with CM. Then  all endomorphisms of $A_f$ are defined over  the Hilbert class field of $K$  and $A_f$ is unique. Let $g_c$ be the dimension of the  abelian variety $A_f$. Due to the fact that $g>1+2 g_c$ ($p\equiv 3 \pmod 4$), resp. $g^+>1+2 g_c$ ($p\equiv 3 \pmod 8)$, the non-existence of an automorphism not defined over $K$ is
  obtained  by applying  the same argument   used in the proof of Lemma 1.4 in  \cite{KM}.

\vskip 0.3 cm

\vskip 0.3 cm
\noindent {\it Step 2: The only primes $N$ which can divide the order of a nontrivial automorphism of $X_{ns}(p)$ or $X_{ns}^+(p)$  are the displayed in the following tables}
\begin{equation}\label{tables}
X_{ns}(p):\quad
\begin{array}{|c|r|}
 p &   N\phantom{cc}\\ \hline
           13 &   2,3, 7\\
           17 &  2,3\\
          19 &   2,3,5\\
         23&  2,3,11\\
         29 & 2,3,5,7\\
          31 &  2,3,5\\
           \hline
          \end{array}\,,\quad \quad X_{ns}^+(p):\quad  \begin{array}{|c|r|}
  p &  N\phantom{cc}\\ \hline
13 &2\\
         17 &  2\\
         19 &  2,3,5\\
         23&  2,3\\
         29 & 2,3,7\\
         31 &  2,3\\\hline
\end{array}
\end{equation}

 The number fields which appear in the decomposition of  $\End_K^0(J_{ns}(p))$ (see (\ref{dc})), resp. $\End_K^0(J_{ns}^+(p))$ (see (\ref{dc+})), only contain the roots of unity $\pm 1$. The only matrix algebras in this decomposition  are of the form $M_2(F)$ for $f\in \cS_{rm}$ and $f\in \cS_{t}$ for $J_{ns}(p)$, resp. $f\in \cS_{rm}^+$ and  $f\in \cS_{t}^+$ for $J_{ns}^+(p)$. In the first case, $F=F_f$ and, in the second case, $F=E_f$. In any case, $F$ is a totally real number field.  If there exists a nontrivial automorphism of order an odd prime $N$, then  the maximal real subfield $K_N$ of the $N$-th cyclotomic field must be contained in some of these number fields $F$. In particular, $N-1$ must divide $ 2 [F:\Q]$. By looking at  the following tables,
 obtained from Table 1,
$$
X_{ns}(p):
\begin{array}{|c|r|}
 p & [F :\Q] \phantom{cc}\\ \hline
           13 & 1,3 \\
           17 &  2,3 \\
          19 &  2,3 \\
         23& 2,5 \\
         29 &1,2,3,6,8 \\
          31 & 1,2,6, 8 \\
           \hline
          \end{array}\,,\quad \quad X_{ns}^+(p): \begin{array}{|c|r|}
  p & [F :\Q] \\ \hline
13 &-\\
         17 &  - \\
         19 &  2 \\
         23& 2 \\
         29 &3 \\
         31 &  8 \\\hline
\end{array}\,,
$$
  we obtain a few possibilities for $N$. After checking all of them, we obtain that the only cases in which  $K_N$ is contained in some $F$ are the displayed in (\ref{tables}).

\vskip 0.3 cm

\noindent {\it Step 3: There are no   automorphisms of $X_{ns}(p)$ and $X_{ns}^+(p)$  of  odd order.}

The claim is obtained applying Theorem \ref{crit} to the curves $X_{ns}(p)\otimes \F_\ell$ and $X_{ns}^+(p)\otimes \F_\ell$, where $\ell$ is a prime splitting in $K$,  and for all $N^m=N$ as in (\ref{tables}). We only show the case $N=3$:
$$\begin{array}{c|c|r|r|}
p  & \ell&X_{ns}(p)\otimes\F_\ell:\sum_nP_3(n) & X_{ns}^+(p)\otimes \F_\ell:\sum_nP_3(n) \\ \hline
13  & 3     & \sum_{n\leq 16}P_3(n)= 12  &      -\\[6 pt]
17  &2 &\sum_{n\leq 34}P_3(n)=18&  -\\ [6 pt]
19  &5 &\sum_{n\leq 31}P_3(n)= 24  & \sum_{n\leq 14}P_3(n)=12\\ [6 pt]
23  &2  &\sum_{n\leq 52}P_3(n)=35 &  \sum_{n\leq 19 }P_3(n)=16\\ [6 pt]
29   &5 &\sum_{n\leq 76}P_3(n)= 58 &  \sum_{n\leq 47}P_3(n)=27\\ [6 pt]
31   &2 &\sum_{n \leq 86}P_3(n)= 66 &  \sum_{n \leq  58 }P_3(n)=31\\
\end{array}
$$

\vskip 0.3 cm
\noindent {\it Step 4: The group  $\Aut(X_{ns}^+(p))$ is trivial.}

We only need to prove that $X_{ns}^+(p)$ does not have any involutions defined over $K$.
Again, the claim is obtained  applying Theorem \ref{crit} to the curves $X=X_{ns}^+(p)\otimes \F_\ell$ for $p\neq 19$ and $X=X_{ns}^+(p)\otimes \F_{\ell^2}$ for $p= 19$,  and  $N^m=2$:
$$\begin{array}{c|c|c|c|}
p & \ell & \sum_n P_2(n)\\ \hline
17 & 2&  \sum_{n\leq 59}P_2(n)=15 \\ [6 pt]
19 & 2 &\sum_{n\leq 83}P_2(n)=19\\ [6 pt]
23 & 2     &\sum_{n\leq  95}P_2(n)=29\\ [6 pt]
29 & 5  &\sum_{n\leq 253}P_2(n)=51 \\ [6 pt]
31 & 2    &\sum_{n\leq 258 }P_2(n)= 59\\
\end{array}
$$
For $p=19$, we have changed the prime $\ell=5$ by   $\ell=2$ ($2$ is inert in $K$), because for $\ell=5$ the sequence $P_2(n)$ turns out to be equal to $0$ for $8\leq n\leq 200$.
\vskip 0.3 cm

\noindent {\it Step 5: The modular involution $w$ is the only nontrivial automorphism  of  $X_{ns}(p)$.}

      A nontrivial automorphism different from $w$   does not commute with $w$ because the group $\Aut(X_{ns}^+(p))$ is trivial. Assume that there is a nontrivial automorphism $u$ of $X_{ns}(p)$ different from $w$. Since the order  of $\Aut_K(X_{ns}(p))$ is a power of $2$,  we can suppose that $u$ is an involution different from  $w$.
The automorphism  $v=u\cdot w$ cannot be an involution, otherwise $u$ and $w$ would commute. Therefore, either $v$ or a power of $v$ has order $4$.

Now, applying Theorem \ref{crit} to $X=X_{ns}(p)\otimes \F_\ell$ for  $N^m=4$ and $p\neq 13, 19$, we obtain
$$\begin{array}{c|c|c|c|}
p & \ell &\sum_nP_4(n)\\ \hline
17 & 2&  \sum_{n\leq 81}P_4(n)=38\\ [6 pt]
%19 & - &--\\ [6 pt]
23 & 2     &\sum_{n\leq 127}P_4(n)=70\\ [6 pt]
29 & 5  &\sum_{n\leq 143}P_4(n)=115\\ [6 pt]
31 & 2    &\sum_{n\leq 291 }P_4(n)= 134\\
\end{array}
$$
Hence, for these four values of $p$ the statement is proved.
 For $p=13$ or  $19$, the sequence $P_4(n)$ turns out to be  equal to  $0$ for $6<n\leq 250$, even changing the prime $\ell$. Nevertheless, applying Theorem \ref{crit} for $N^m=8$, we prove that $X_{ns}(13)$ an $X_{ns}(19)$  do not  have any automorphisms of order $8$:
 $$\begin{array}{c|c|c|c|}
p & \ell & \sum_n P_8(n)\\ \hline
13 & 3&  \sum_{n\leq 15}P_8(n)= 34\\ [6 pt]
19 & 5 &\sum_{n\leq 34}P_8(n)=58
\end{array}
$$
Therefore, the order of any automorphism of $X_{ns}(p)$ must divide $4$. Assume that there is  $v\in\Aut(X_{ns}(p))$ of order $4$. Then, the automorphism $u:=v^2\cdot w$ can only have order $2$ or $4$. On the one hand,  $u$ cannot be an involution since  $v^2$ and $w$ do no commute. On the other hand, if $u$ has order $4$, then $u^2$ is an involution different from $w$ and, thus, $u^2\cdot w=v^2\cdot w\cdot v^2$  must have order $4$, but $(u^2\cdot w)^2=1$. Therefore, none of these two curves has automorphisms of order $4$.
\hfill $\Box$

\bibliography{X_ns(11)}{}
\bibliographystyle{alpha}

\vskip 0.4 cm

%\begin{footnotesize}
\begin{tabular}{l}
Josep Gonz\'alez\\
\texttt{josepg\,\footnotesize{$@$}\,ma4.upc.edu}\\
Departament de Matem\`atica Aplicada 4  \\
Universitat Polit\`ecnica de Catalunya  \\
EPSEVG, Avinguda V\'ictor Balaguer 1\\
08800 Vilanova i la Geltr\'u, Spain\\[20pt]
\end{tabular}

\end{document}